\documentclass{svmult}

\usepackage{mathptmx}        
\usepackage{helvet}          
\usepackage{courier}         

\usepackage{makeidx}         
\usepackage{graphicx}        
\usepackage{multicol}        
\usepackage[bottom]{footmisc}
\usepackage{mathtools}
\usepackage{xcolor}

\makeindex             

\usepackage{subeqnarray}
\usepackage{amsmath}
\usepackage{amssymb}

\begin{document}

\title*{Almost Diagonalization of Pseudodifferential Operators}
\author{S. Ivan Trapasso}
\institute{
Salvatore Ivan Trapasso
\at {Dipartimento di Scienze Matematiche ``G. L. Lagrange", Politecnico di Torino \\ Corso Duca degli Abruzzi 24, 10129 Torino (Italy)} \\
\email{salvatore.trapasso@polito.it}
}

\maketitle
\abstract{In this review we focus on the almost diagonalization of pseudodifferential operators and highlight the advantages that time-frequency techniques provide here. In particular, we retrace the steps of an insightful paper by Gr\"ochenig, who succeeded in characterizing a class of symbols previously investigated by Se\"ostrand by noticing that Gabor frames almost diagonalize the corresponding Weyl operators. This approach also allows to give new and more natural proofs of related results such as boundedness of operators or algebra and Wiener properties of the symbol class. Then, we discuss some recent developments on the theme, namely an extension of these results to a more general family of pseudodifferential operators and similar outcomes for a symbol class closely related to Sj\"ostrand's one.}
\vfill
This is a pre‐copyedited version of a contribution published in \textit{Landscapes of Time-Frequency Analysis} (Boggiatto P. et al. (eds)) published by Birkh\"auser, Cham. The definitive authenticated version is available online via https://doi.org/10.1007/978-3-030-05210-2\_14.

\vfill


\pagebreak
\section{Introduction}

 The wide range of problems that one can tackle by means of Time-frequency Analysis bears witness to the relevance of this quite modern discipline stemmed from both pure and applied issues in harmonic analysis. There is no way to provide here a comprehensive bibliography on the theme, which would encompass studies in quantum mechanics and partial differential equations. We confine ourselves to list some references to be used as points of departure for a walk through the topic: see \cite{BogetalTRANS,CNStricharzJDE2008,deGossonDiasPrata2014,deGossonWigner2017,Ruzhansky2016,SugimotoWang2011}. Besides the countless achievements as tool for other fields, Gabor analysis is a fascinating subject in itself and it may happen to shed new light on established facts in an effort to investigate the subtle problems underlying its foundation. We report here the case of Gr\"{o}chenig's work \cite{Grochenig_2006_Time}: the author retrieved and extended well-known outcomes obtained by Sj\"{o}strand within the realm of ``hard" analysis - cf. \cite{Sjo94,Sjo95}, and this was achieved using techniques from phase space analysis. We will give a detailed account in the subsequent sections, but let us briefly introduce here the main characters of this story. 
 \\ 
The (cross-)Wigner distribution is a quadratic time-frequency representation of signals $f,g$ in suitable function spaces (for instance $f,g\in\mathcal{S}(\mathbb{R}^{d})$, the Schwartz class) defined as 
\begin{equation}
W(f,g)(x,\omega )=\int_{\mathbb{R}^{d}}e^{-2\pi iy \omega }f\left(x+\frac{y}{2}\right)%
\overline{g\left(x-\frac{y}{2}\right)}\,\D y.
\label{wig}
\end{equation}
It is possible to associate a pseudodifferential operator to this representation, namely the so-called Weyl transform - it is a quite popular quantization rule in Physics community. Given a tempered distribution $\sigma\in\mathcal{S}'(\mathbb{R}^{2d})$ as \emph{symbol} (also \emph{observable}, in physics vocabulary), the corresponding Weyl transform maps $\mathcal{S}(\mathbb{R}^d)$ into $\mathcal{S}'(\mathbb{R}^{d})$ and can be defined via duality by
\begin{equation}
\langle \mathrm{Op_W} (\sigma)f,g\rangle =
\langle
\sigma,W(g,f)\rangle, \quad f,g\in \mathcal{S}(\mathbb{R}^{d}).
\label{weylt}
\end{equation}

The Weyl transform has been thoroughly studied in \cite{Grochenig_2001_Foundations,WongWeylTransform1998} among others. 
In his aforementioned works, Sj\"ostrand proved that Weyl operators with symbols of special type satisfy a number of interesting properties concerning their boundedness and algebraic structure as a set. In terms that will be specified later in Section 3, we can state that the set of such operators is a spectral invariant *-subalgebra of $\mathcal{B}(L^2(\mathbb{R}^{d}))$, the ($C*$-)algebra of bounded operators on $L^2(\mathbb{R}^{d})$. 

To be precise, given a Schwartz  function $g\in\mathcal{S}(\mathbb{R}^{2d})\setminus \{0\}$, we provisionally define the \emph{Sj\"ostrand's class} as the space of tempered distributions $\sigma\in\mathcal{S}'(\mathbb{R}^{2d})$ such that
$$ \int_{\mathbb{R}^{2d}}\sup_{z\in\mathbb{R}^{2d}} |\langle \sigma, \pi(z,\zeta)g\rangle| \D \zeta<\infty.
$$
As a rule of thumb, notice that a symbol in $M^{\infty,1}(\mathbb{R}^{2d})$ locally (i.e. for fixed $z\in \mathbb{R}^{2d}$) coincides with the Fourier transform of a $L^1$ function. Furthermore, it can be proved that this somewhat exotic symbol class contains classical H\"ormander's symbols of type $S^0_{0,0}$, together with non-smooth ones. 

The crucial remark here is that Sj\"ostrand's class actually coincides with a function space of a particular type, namely the modulation space $M^{\infty,1}(\mathbb{R}^{2d})$. In more general terms, modulation spaces (and also related Wiener amalgam spaces of special type) are Banach spaces defined by means of estimates on time-frequency concentration and decay of its elements - see Section 2 for the details. They were introduced by Feichtinger
in the '80s (cf. the pioneering papers \cite{Segal81.Feichtinger_1981_Banach,feichtinger1983modulation}) and soon established themselves as the optimal environment for time-frequency analysis. Nevertheless, they also provide a fruitful context to set problems in harmonic analysis and PDEs - see for instance \cite{deGossonsymplectic2011,deGossonGRomero2016,Wangbook2011}. 

Gr\"ochenig deeply exploited this connection with time-frequency analysis by proving that Sj\"ostrand's results extend to more general modulation spaces and, more importantly, he was able to completely characterize symbols in these classes by means of a property satisfied by the corresponding Weyl operators, namely approximate diagonalization. This is a classical problem in pure and applied harmonic analysis - a short list of references is \cite{cordero2013wiener,Labate2008, meyer1990ondelettes,rochberg1998pseudodifferential}. 
We will thoroughly examine Gr\"ochenig's results in Section 3. Here, we limit ourselves to heuristically argue that the choice of a certain type of symbols assures that the corresponding Weyl operators preserve the time-frequency localization, since their ``kernel" with respect to continuous or discrete time-frequency shifts satisfies a convenient decay condition. 

In the subsequent Section 4 we report some results on almost diagonalization obtained by the author in a recent joint work with Elena Cordero and Fabio Nicola - see \cite{CNT18}.
Mimicking the scheme which leads to define the Weyl transform, in \cite{BogetalTRANS} the authors consider a one-parameter family of time-frequency representations (\emph{$\tau$-Wigner distributions}) and  also define the corresponding pseudodifferential operators $\mathrm{Op}_{\tau}$ via duality. Precisely, for $\tau \in \lbrack 0,1]$, the (cross-)$\tau $-Wigner distribution is given by
\begin{equation}
W_{\tau }(f,g)(x,\omega )=\int_{\mathbb{R}^{d}}e^{-2\pi iy\zeta }f(x+\tau y)%
\overline{g(x-(1-\tau )y)}\,\D y,\quad f,g\in \mathcal{S}(\mathbb{R}^{d}),
\label{tauwig}
\end{equation}
whereas the corresponding $\tau$-pseudodifferential operator is  defined by
\begin{equation}
\langle \mathrm{Op}_{\tau}  (a)f,g\rangle =
\langle
a,W_{\tau }(g,f)\rangle, \quad f,g\in \mathcal{S}(\mathbb{R}^{d}).
\label{tauweak}
\end{equation}
For $\tau=1/2$ we recapture the Weyl transform and the usual Wigner distribution, while the cases $\tau=0,1$ respectively cover the classical theory of Kohn-Nirenberg and anti-Kohn-Nirenberg operators - whose corresponding distributions are also known as Rihaczek and conjugate-Rihaczek distributions respectively.

Our contribution aims at enlarging the area of application of Gr\"ochenig's result along two directions. First, one finds that symbols in the Sj\"ostrand's class are in fact characterized by almost diagonalization of the corresponding $\tau$-pseudodifferential operators for any $\tau \in [0,1]$. While this is not surprising for reasons that will be discussed later, it seems worthy of interest to get a similar result for symbols belonging to a function space closely related to $M^{\infty,1}$, namely the Wiener amalgam space $W(\mathcal{F} L^{\infty},L^1)$. The connection between these spaces is established by Fourier transform: in fact, the latter exactly contains the Fourier transforms of symbols in the Sj\"ostrand's class. It is important to remark that even if the spirit of the result is the same, numerous differences occur and we try to clarify the intuition behind this situation in Section 4.

To conclude, we take advantage of this characterization in regards to boundedness results. We were able to study the boundedness of $\tau$-pseudodifferential operator covering several possible choices among modulation and Wiener amalgam space for symbols classes and spaces on which they act. We mention that in a number of these outcomes we have benefited from a strong linkage with the theory of Fourier integral operators. Besides, the latter condition also made possible to establish (or disprove) the algebraic properties considered by Sj\"ostrand for special classes of $\tau$-pseudodifferential operators.

\section{Preliminaries}
\noindent
 \textbf{Notation.} We write $t^2=t\cdot t$, for $t\in\mathbb{R}^{d}$, and
$xy=x\cdot y$ is the scalar product on $\mathbb{R}^{d}$. The Schwartz class is denoted by  $\mathcal{S}(\mathbb{R}^{d})$, the space of tempered
distributions by  $\mathcal{S}'(\mathbb{R}^{d})$.  The brackets  $\langle
f,g\rangle$ denote both the duality pairing between $\mathcal{S}' (\mathbb{R}^{d})$ and $\mathcal{S} (\mathbb{R}^{d})$ and the inner product $\langle f,g\rangle=\int f(t){\overline {g(t)}} \D t$ on $L^2(\mathbb{R}^{d})$. In particular, we assume it to be conjugate-linear in the second argument.
The symbol $\lesssim$ means that the underlying inequality holds up to a positive constant factor $C>0$ on the RHS:
$$ f\lesssim g\quad\Rightarrow\quad\exists C>0\,:\,f\le Cg. $$.
The Fourier transform of a function $f$ on $\mathbb{R}^{d}$ is normalized as
\[
\mathcal{F} f(\xi)= \int_{\mathbb{R}^{d}} e^{-2\pi i x\xi} f(x)\, \D x.
\]

Given $x,\omega \in \mathbb{R}^{d}$, the modulation $M_{\omega}$ and translation $T_{x}$ operators acts on a function $f$ (on $\mathbb{R}^{d}$) as 
\[
M_{\omega}f\left(t\right)= e^{2\pi it \omega}f\left(t\right),\qquad T_{x}f\left(t\right)= f\left(t-x\right).
\]
We write a point in phase space as
$z=(x,\omega)\in\mathbb{R}^{2d}$, and  the corresponding phase-space shift 
acting on a function or distribution as
\begin{equation}
\label{eq:kh25}
\pi (z)f(t) = e^{2\pi i \omega t} f(t-x), \, \quad t\in\mathbb{R}^{d}.
\end{equation}

Denote by $J$ the canonical symplectic matrix in $\mathbb{R}^{2d}$:
\[
J=\left(\begin{array}{cc}
0_{d\times d} & I_{d\times d}\\
-I_{d\times d} & 0_{d\times d}
\end{array}\right)\in\mathrm{Sp}\left(2d,\mathbb{R}\right),
\]
where  the
symplectic group $\mathrm{Sp}\left(2d,\mathbb{R}\right)$ is defined
by
$$
\mathrm{Sp}\left(2d,\mathbb{R}\right)=\left\{M\in \mathrm{GL}(2d,\mathbb{R}):\;M^{\top}JM=J\right\}.
$$
Observe that, for $z=\left(z_{1},z_{2}\right)\in\mathbb{R}^{2d}$, we have
$Jz=J\left(z_{1},z_{2}\right)=\left(z_{2},-z_{1}\right),$  $J^{-1}z=J^{-1}\left(z_{1},z_{2}\right)=\left(-z_{2},z_{1}\right)=-Jz,$ and 
$J^{2}=-I_{2d\times2d}$.
\\ \\
\textbf{Short-time Fourier transform.} Let $f\in\mathcal{S}'(\mathbb{R}^{d})$ and $g\in \mathcal{S}(\mathbb{R}^{d}) \setminus \{0\}$. The short-time Fourier transform (STFT) of $f$ with window function $g$ is defined as
\begin{equation}\label{STFTdef}
V_gf(x,\omega )=\langle f,\pi(x,\omega ) g\rangle=\mathcal{F} (fT_x g)(\omega)=\int_{\mathbb{R}^{d}}
f(y)\, {\overline {g(y-x)}} \, e^{-2\pi iy \omega }\, \D y.
\end{equation}
We remark that the last expression has to be intended in formal sense, but it truly represents the integral corresponding to the inner product $\langle f,\pi(x,\omega ) g\rangle$ whenever $f,g\in L^2(\mathbb{R}^{d})$. 

Recall the fundamental property of time-frequency analysis:
\begin{equation}\label{FI}
V_{g}f\left(x,\omega\right)=e^{-2\pi ix\omega}V_{\hat{g}}\hat{f}\left(J(x,\omega)\right).
\end{equation}

\textbf{Gabor frames. }Let $\Lambda=A\mathbb{Z}^{2d}$,  with $A\in \mathrm{GL}(2d,\mathbb{R})$, be a lattice in the time-frequency plane.
The set  of time-frequency shifts $\mathcal{G}(\varphi,\Lambda)=\{\pi(\lambda)\varphi:\
\lambda\in\Lambda\}$ for a  non-zero $\varphi\in L^2(\mathbb{R}^{d})$ (the so-called window function) is called Gabor system. A Gabor system $\mathcal{G}(\varphi,\Lambda)$ is said 
to be a Gabor frame if the lattice is such thick that the energy content of a signal as sampled on the lattice by means of STFT is comparable with its total energy, that is: there exist 
constants $A,B>0$ such that
\begin{equation}\label{gaborframe}
A\|f\|_2^2\leq\sum_{\lambda\in\Lambda}|\langle f,\pi(\lambda)\varphi\rangle|^2\leq B\|f\|^2_2,\qquad \forall f\in L^2(\mathbb{R}^{d}).
\end{equation}

\subsection{Function spaces} 
\noindent
\textbf{Weight functions.} Let us call \emph{admissible
weight function} any non-negative continuous
function $v$ on $\mathbb{R}^{2d}$ such that: 
\begin{enumerate}
	\item $v\left(0\right)=1$ and $v$ is even in each coordinate: 
	\[
	v\left(\pm z_{1},\ldots,\pm z_{2d}\right)=v\left(z_{1},\ldots,z_{2d}\right).
	\]
	\item $v$ is submultiplicative, that is 
	\[
	v\left(w+z\right)\le v\left(w\right)v\left(z\right)\qquad\forall w,z\in\mathbb{R}^{2d}.
	\]
	\item $v$ satisfies the Gelfand-Raikov-Shilov (GRS) condition:
	\begin{equation}\label{GRS}
	\lim_{n\rightarrow\infty}v\left(nz\right)^{\frac{1}{n}}=1\qquad\forall z\in\mathbb{R}^{2d}.
	\end{equation}
\end{enumerate}

Examples of admissible weights are given by $v\left(z\right)=e^{a\left|z\right|^{b}}\left(1+\left|z\right|\right)^{s}\log^{r}\left(e+\left|z\right|\right)$,
with real parameters $a,r,s\ge0$ and $0\leq b<1$. Functions of polynomial growth such as
\begin{equation}\label{vs}
v_{s}\left(z\right)=\left\langle z\right\rangle ^{s}=\left(1+\left|z\right|^{2}\right)^{\frac{s}{2}},\qquad z\in\mathbb{R}^{2d},\,s\ge0
\end{equation}
are admissible weights too. From now on, $v$ will denote an admissible weight function unless otherwise specified. 
We remark that the GRS condition is exactly the technical tool required to forbid an exponential growth of the weight in some direction. For further discussion on this feature, see \cite{massop}. \\
Given a submultiplicative weight $v$, a positive function $m$ on $\mathbb{R}^{2d}$ is called \emph{$v$-moderate weight} if there exists a constant $C\geq 0$ such that
     $$
         m(z_1+z_2)\leq Cv(z_1)m(z_2)\;,\qquad z_1, z_2\in\mathbb{R}^{2d}.
     $$
The set of all $v$-moderate weights will be denoted by $\mathcal{M}_v(\mathbb{R}^{2d})$.
     
In order to remain in the framework of tempered distributions, in what follows we shall always assume that weight functions $m$ on $\mathbb{R}^{d}$ under our consideration satisfy the following condition:
\begin{equation}\label{M}
m(z)\geq 1,\quad \forall z\in \mathbb{R}^{d}\quad \mbox{or} \quad m(z)\gtrsim \langle z\rangle^{-N},
\end{equation}
for some $N\in\mathbb{N}$. The same holds with suitable modifications for weights on $\mathbb{R}^{2d}$.
\\ \\ 
\noindent
\textbf{Modulation spaces. }
Given a non-zero window $g\in\mathcal{S}(\mathbb{R}^{d})$, a $v$-moderate weight
function $m$ on $\mathbb{R}^{2d}$ satisfying \eqref{M}, and $1\leq p,q\leq
\infty$, the {\it
	modulation space} $M^{p,q}_m(\mathbb{R}^{d})$ consists of all tempered
distributions $f\in\mathcal{S}'(\mathbb{R}^{d})$ such that $V_gf\in L^{p,q}_m(\mathbb{R}^{2d} )$
(weighted mixed-norm space). The norm on $M^{p,q}_m$ is
$$
\|f\|_{M^{p,q}_m}=\|V_gf\|_{L^{p,q}_m}=\left(\int_{\mathbb{R}^{d}}
\left(\int_{\mathbb{R}^{d}}|V_gf(x,\omega)|^pm(x,\omega)^p\,
\D x\right)^{q/p} \D \omega\right)^{1/q}  \, ,
$$
with suitable modifications if $p=\infty$ or $q=\infty$.
If $p=q$, we write $M^p_m$ instead of $M^{p,p}_m$, and if $m(z)\equiv 1$ on $\mathbb{R}^{2d}$, then we write $M^{p,q}$ and $M^p$ for $M^{p,q}_m$ and $M^{p,p}_m$.

It can be proved (see \cite{Grochenig_2001_Foundations}) that $M^{p,q}_m (\mathbb{R}^{d} )$ is a Banach space whose definition is independent of the choice of the window $g$ - meaning that different windows provide equivalent norms on $M^{p,q}_m$. The window class can be extended to $M^1_v$, cf. \cite[Thm.~11.3.7]{Grochenig_2001_Foundations}. Hence, given any $g \in M^1_v
(\mathbb{R}^{d} )$ and $f\in M^{p,q}_m $ we have 
\begin{equation}\label{normwind}
\|f\|_{M^{p,q}_m } \asymp \|V_{g}f \|_{L^{p,q}_m }.
\end{equation}
We  recall the inversion formula for
the STFT (see  \cite[Proposition 11.3.2]{Grochenig_2001_Foundations}). If $g\in M^{1}_v(\mathbb{R}^{d})\setminus\{0\}$,
$f\in M^{p,q}_m(\mathbb{R}^{d})$, with $m$ satisfying \eqref{M}, then
\begin{equation}\label{invformula}
f=\frac1{\|g\|_2^2}\int_{\mathbb{R}^{2d}} V_g f(z) \pi (z)  g\, \D z \, ,
\end{equation}
and the  equality holds in $M^{p,q}_m(\mathbb{R}^{d})$.\par
The adjoint operator of $V_g$,  defined by
$$V_g^\ast F(t)=\int_{\mathbb{R}^{2d}} F(z)  \pi (z) g \D z \, ,
$$
maps the Banach space $L^{p,q}_m(\mathbb{R}^{2d})$ into $M^{p,q}_m(\mathbb{R}^{d})$. In particular, if $F=V_g f$ the inversion formula \eqref{invformula} reads
\begin{equation}\label{treduetre}
{\rm Id}_{M^{p,q}_m}=\frac 1 {\|g\|_2^2} V_g^\ast V_g.
\end{equation}

\noindent
\textbf{Wiener Amalgam Spaces.} Fix $g\in \mathcal{S}(\mathbb{R}^{d}) \setminus \left\{ 0 \right\} $ and consider \emph{even} weight functions $u,w$ on $\mathbb{R}^{d}$ satisfying \eqref{M}. The Wiener amalgam space $W(\mathcal{F} L^p_u,L^q_w)(\mathbb{R}^{d})$ is the space of distributions $f\in\mathcal{S}'(\mathbb{R}^{d})$ such that
\[
\|f\|_{W(\mathcal{F} L^p_u,L^q_w)(\mathbb{R}^{d})}:=\left(\int_{\mathbb{R}^{d}}
\left(\int_{\mathbb{R}^{d}}|V_gf(x,\omega)|^p u^p(\omega)\,
\D \omega \right)^{q/p} w^q(x)\D x\right)^{1/q}<\infty  \,
\]
with obvious modifications for $p=\infty$ or $q=\infty$.

Using the fundamental identity of time-frequency analysis\, \eqref{FI}, we have $|V_g f(x,\omega)|=|V_{\hat g} \hat f(\omega,-x)| = |\mathcal F (\hat f \, T_\omega \overline{\hat g}) (-x)|$  and (since $u(x)=u(-x)$)
$$
\| f \|_{{M}^{p,q}_{u\otimes w}} = \left( \int_{\mathbb{R}^{d}} \| \hat f \ T_{\omega} \overline{\hat g} \|_{\mathcal{F} L^p_u}^q w^q(\omega) \ \D \omega \right)^{1/q}
= \| \hat f \|_{W(\mathcal{F} L_u^p,L_w^q)}.
$$

Hence the Wiener amalgam spaces under our consideration are simply the image under Fourier transform\, of modulation spaces with weights of tensor product type, namely $m(x,\omega)=u\otimes w (x,\omega)=u(x)w(\omega)$:
\begin{equation}\label{W-M}
\mathcal{F} ({M}^{p,q}_{u\otimes w})=W(\mathcal{F} L_u^p,L_w^q).
\end{equation}

For this reason among others, their inventor H. Feichtinger suggested to call them modulation spaces too - although in a generalized sense, see \cite{Feich2006} for an intriguing conceptual account on the theme.

\subsection{$\tau$-Pseudodifferential Operators}
Let us introduce the $\tau$-pseudodifferential operators as it is customary in time-frequency analysis, i.e. by means of superposition of time-frequency shifts:
\begin{equation}
\mathrm{Op}_{\tau}  \left(\sigma\right)f\left(x\right) = \int_{\mathbb{R}^{2d}}\hat{\sigma}\left(\omega,u\right)e^{-2 \pi i \left(1-\tau\right)\omega u}\left(T_{-u}M_{\omega}f\right)\left(x\right)\D u \D \omega,\qquad x\in\mathbb{R}^{d},\label{opt tfs}
\end{equation}
for any $\tau\in\left[0,1\right]$. The symbol $\sigma$ and the function
$f$ belong to suitable function spaces, to be determined in order
for the previous expression to make sense. As an example, minor modifications to \cite[Lem. 14.3.1]{Grochenig_2001_Foundations} give that $\mathrm{Op}_{\tau}  \left(\sigma\right)$
maps $\mathcal{S}\left(\mathbb{R}^{d}\right)$ to $\mathcal{S}'\left(\mathbb{R}^{2d}\right)$
whenever $\sigma\in\mathcal{S}'\left(\mathbb{R}^{2d}\right)$. 

Assuming that \eqref{opt tfs} is a well-defined absolutely convergent integral (for instance, it is enough to assume $\hat{\sigma}\in L^{1}\left(\mathbb{R}^{2d}\right)$), easy computations lead to the usual integral form of $\tau$-pseudodifferential
operators, namely
\begin{flalign*}
\mathrm{Op}_{\tau}  \left(\sigma\right)f\left(x\right) & =\int_{\mathbb{R}^{2d}}e^{2\pi i\left(x-y\right)\omega}\sigma\left(\left(1-\tau\right)x+\tau y,\omega\right)f\left(y\right)\D y \D \omega.
\end{flalign*}

We finally aim to represent $\mathrm{Op}_{\tau} (\sigma)$ as an integral operator
of the form 
\[
\mathrm{Op}_{\tau}  \left(\sigma\right)f\left(x\right)=\int_{\mathbb{R}^{2d}}k\left(x,y\right)f\left(y\right)\D y.
\]
Let us introduce the operator $\mathfrak{T}_{\tau}$ acting on functions
on $\mathbb{R}^{2d}$ as
\[
\mathfrak{T}_{\tau}F\left(x,y\right) = F\left(x+\tau y,x-\left(1-\tau\right)y\right),\qquad\mathfrak{T}_{\tau}^{-1}F\left(x,y\right)=F\left(\left(1-\tau\right)x+\tau y,x-y\right),
\]
and denote by $\mathcal{F}_{i}$, $i=1,2$, the partial Fourier
transform with respect to the $i$-th $d-$dimensional variable (it is then clear that $\mathcal{F}=\mathcal{F}_{1}\mathcal{F}_{2}$).

Since the operators $\mathfrak{T}_{\tau}$ and $\mathcal{F}_{i}$
are continuous bijections on $\mathcal{S}\left(\mathbb{R}^{2d}\right)$,
the kernel $k$ is well-defined (as a tempered distribution) also
for symbols in $\mathcal{S}'\left(\mathbb{R}^{2d}\right)$ and we
finally recover the representation by duality given in the Introduction according to \cite{BogetalTRANS}.

\begin{proposition}{\label{optker}}
	For any symbol $\sigma\in\mathcal{S}'\left(\mathbb{R}^{2d}\right)$
	and any real $\tau\in\left[0,1\right]$, the map $\mathrm{Op}_{\tau} (\sigma) \,:\,\mathcal{S}\left(\mathbb{R}^{d}\right)\rightarrow\mathcal{S}\left(\mathbb{R}^{d}\right)$
	is defined as integral operator with distributional kernel 
	\[
	k=\mathfrak{T}_{\tau}^{-1}\mathcal{F}_{2}^{-1}\sigma\in\mathcal{S}'\left(\mathbb{R}^{2d}\right),
	\]
	meaning that, for any $f,g\in\mathcal{S}\left(\mathbb{R}^{d}\right)$,
	\[
	\left\langle \mathrm{Op}_{\tau}  \left(\sigma\right)f,g\right\rangle =\left\langle k,g\otimes\overline{f}\right\rangle .
	\]
	In particular, since the representation 
	\[
	W_{\tau}\left(f,g\right)\left(x,\omega\right)=\mathcal{F}_{2}\mathfrak{T}_{\tau}\left(f\otimes\overline{g}\right)\left(x,\omega\right)
	\]
	holds for $f,g\in\mathcal{S}\left(\mathbb{R}^{d}\right)$, we have
	\[
	\left\langle \mathrm{Op}_{\tau}  \left(\sigma\right)f,g\right\rangle =\left\langle \sigma,W_{\tau}\left(g,f\right)\right\rangle .
	\]
\end{proposition}

As a consequence of the celebrated Schwartz's kernel theorem (see for instance \cite[Theorem 14.3.4]{Grochenig_2001_Foundations}), we are able to relate the representations for $\tau$-pseudodifferential operators given insofar. 

\begin{theorem}\label{optreps}
	Let $T\,:\,\mathcal{S}\left(\mathbb{R}^{d}\right)\rightarrow\mathcal{S}'\left(\mathbb{R}^{d}\right)$
	be a continuous linear operator. There exist tempered distributions
	$k,\sigma,F\in\mathcal{S}'\left(\mathbb{R}^{d}\right)$ and $\tau\in\left[0,1\right]$
	such that $T$ admits the following representations:
	\begin{enumerate}
		\item[$(i)$] as an integral operator: $\left\langle Tf,g\right\rangle =\left\langle k,g\otimes\overline{f}\right\rangle $
		for any $f,g\in\mathcal{S}\left(\mathbb{R}^{d}\right)$;
		\item[$(ii)$] as a $\tau$-pseudodifferential operator $T=\mathrm{Op}_{\tau}\left(\sigma\right)$
		with symbol $\sigma$;
		\item[$(iii)$] as a superposition (in a weak sense) of time-frequency shifts : $$T=\int_{\mathbb{R}^{2d}}F\left(x,\omega\right)e^{2\left(1-\tau\right)\pi ix\omega}T_{x}M_{\omega}\D x \D \omega.$$
	\end{enumerate}
	The relations among $k$, $\sigma$ and $F$ are the following:
	\[
	\sigma=\mathcal{F}_{2}\mathfrak{T}_{\tau}k,\qquad F=\mathcal{I}_{2}\hat{\sigma},
	\]
	where $\mathcal{I}_{2}$ denotes the reflection in the second $d$-dimensional
	variable (i.e. $\mathcal{I}_{2}G\left(x,\omega\right)=G\left(x,-\omega\right)$,
	$\left(x,\omega\right)\in\mathbb{R}^{2d}$). 
\end{theorem}

To conclude this anthology, since the algebraic properties of pseudodifferential operators families will be considered, recall that the composition of Weyl transforms provides a bilinear form on symbols, the so-called \emph{twisted product}:
$$\mathrm{Op_W}(\sigma) \circ \mathrm{Op_W}(\rho) = \mathrm{Op_W} (\sigma \sharp \rho).$$
Although explicit formulas for the twisted product of symbols can be derived (cf. \cite{WongWeylTransform1998}), we will not need them hereafter. Anyway, this is a fundamental notion in order to establish an algebra structure on symbol spaces: it is quite natural to ask if the composition of operators with symbols in the same class reveals to be an operator of the same type for some symbol in the same class. 
Also recall that taking the adjoint of a Weyl operator provides an involution on the level of symbols, since  $\left(\mathrm{Op_W}(\sigma)\right)^*=\mathrm{Op_W}(\overline{\sigma})$.

\section{Time-frequency analysis of the Sj\"ostrand's class}

The study of pseudodifferential operators has a wide and long tradition in the field of mathematical analysis, starting from the monumental work of H\"ormander. It has to be noticed that the classical symbol classes considered in these investigations are usually defined by means of differentiability conditions. In the spirit of time-frequency analysis, we hereby employ modulation and Wiener amalgam spaces as reservoirs of symbols for pseudodifferential operator and hence the short-time Fourier transform to shape the desired properties. \\
Recall that the Sj\"ostrand's class is the modulation
space $M^{\infty,1}(\mathbb{R}^{2d})$ consisting of distributions $\sigma\in\mathcal{S}'(\mathbb{R}^{2d})$ such that
$$ \int_{\mathbb{R}^{2d}} \sup_{z\in\mathbb{R}^{2d}} |\langle \sigma, \pi(z,\zeta)g\rangle| \D \zeta<\infty.
$$

The control on symbols can be improved by weighting the condition on their short-time Fourier transform, i.e. the modulation space norm. In the following we will employ weight functions of type $1\otimes v$, where $v$ is an admissible weight on $\mathbb{R}^{2d}$, according to the properties assumed in the Preliminaries. Weighted Sj\"ostrand's classes of this type are thus defined as

$$M^{\infty,1}_{1\otimes v}\left(\mathbb{R}^{2d}\right)=\left\{ \sigma \in \mathcal{S}' \left(\mathbb{R}^{2d}\right) \, : \, \int_{\mathbb{R}^{2d}}\sup_{z\in \mathbb{R}^{2d}} |V_g\sigma(z,\zeta)|v\left(\zeta\right) \D \zeta < \infty\right\}.$$

A function space closely related to the previous one is the Wiener amalgam space $W(\mathcal{F} L^\infty, L^1_v)(\mathbb{R}^{2d})$. As discussed in the previous section, we have indeed $W(\mathcal{F} L^\infty, L^1_v)(\mathbb{R}^{2d})=\mathcal{F} M^{\infty,1}_{1\otimes v}(\mathbb{R}^{2d})$. Heuristically, a symbol in $W(\mathcal{F} L^\infty, L^1)(\mathbb{R}^{2d})$ locally coincides with the Fourier transform of a $L^\infty(\mathbb{R}^{2d})$ signal and exhibits global decay of $L^1$ type. For instance, the $\delta$ distribution  (in $\mathcal{S}'(\mathbb{R}^{2d})$) belongs to $W(\mathcal{F} L^\infty, L^1)(\mathbb{R}^{2d})$.

Although Sj\"ostrand's definition of the eponym symbol class was quite different from the one given here in terms of modulation spaces, in his works \cite{Sjo94,Sjo95} he proved three fundamental results on Weyl operators with symbols in $M^{\infty,1}$. 

\begin{theorem} ~ 
	\begin{enumerate}
		\item[$(i)$] (\textbf{Boundedness}) If $\sigma\in M^{\infty,1}\left(\mathbb{R}^{2d}\right)$, then $\mathrm{Op_W}(\sigma)$ is a bounded operator on $L^2(\mathbb{R}^{d})$.  
		\item[$(ii)$] (\textbf{Algebra property}) If $\sigma_1, \sigma_2 \in M^{\infty,1}\left(\mathbb{R}^{2d}\right)$ and $\mathrm{Op_W}(\rho)=\mathrm{Op_W}(\sigma_1)\mathrm{Op_W}(\sigma_2)$, then $\rho=\sigma_1 \sharp \sigma_2\in M^{\infty,1}\left(\mathbb{R}^{2d}\right)$.
		\item[$(iii)$] (\textbf{Wiener property}) If $\sigma\in M^{\infty,1}\left(\mathbb{R}^{2d}\right)$ and $\mathrm{Op_W}(\sigma)$ is invertible on $L^2(\mathbb{R}^{d})$,  then $\left[\mathrm{Op_W}(\sigma)\right]^{-1}=\mathrm{Op_W}(\rho)$ for some $\rho\in M^{\infty,1}\left(\mathbb{R}^{2d}\right)$.
	\end{enumerate}
\end{theorem}

For sake of conciseness, we can resume the preceding outcomes by saying that the family of Weyl operators with symbols in Sj\"ostrand's class (denoted by $\mathrm{Op_W}(M^{\infty,1})$) is an inverse-closed Banach *-subalgebra of $\mathcal{B} (L^2(\mathbb{R}^{d}))$. 

Both these results and their original proofs might appear fairly technical at first glance. Nonetheless, they unravel a deep and fascinating analogy between Weyl operators with symbols in the Sj\"ostrand's class and Fourier series with $\ell^1$ coefficients. Similarities of this kind come under the multifaceted problem of spectral invariance, a topic thoroughly explored by Gr\"ochenig in his insightful lecture \cite{massop}.

In view of the structure of $\tau$-pseudodifferential operators as superposition of time-frequency shifts (cf. Equation \eqref{opt tfs}), it can be fruitful to study how operators interact with time-frequency shifts. A measure of this interplay is given by the entries of the infinite matrix which we are going to refer to  as \emph{channel matrix}, according to traditional nomenclature in applied contexts like data transmission. First, fix a non-zero window $\varphi\in M^1_v(\mathbb{R}^{d})\left(\mathbb{R}^{d}\right)$ and a lattice $\Lambda = A\mathbb{Z}^{2d}\subseteq \mathbb{R}^{2d}$, where $A\in\mathrm{GL}(2d,\mathbb{R}) $, such that 
$\mathcal{G}\left(\varphi,\Lambda\right)$ is a Gabor frame  for $L^{2}\left(\mathbb{R}^{d}\right)$. Therefore, the entries of the channel matrix are given by
 $$\langle \mathrm{Op_W}(\sigma)\pi(z)\varphi,\pi(w)\varphi\rangle,\qquad z,w\in\mathbb{R}^{2d}, $$ or $$M(\sigma)_{\lambda,\mu}\coloneqq \langle \mathrm{Op_W}(\sigma)\pi(\lambda)\varphi,\pi(\mu)\varphi\rangle,\qquad \lambda,\mu\in\Lambda,$$ if we restrict to the lattice $\Lambda$. 
 In this context, we could say that $\mathrm{Op_W}$ is almost diagonalized by the Gabor frame $\mathcal{G}(\varphi,\Lambda)$ if its channel matrix exhibits a suitable off-diagonal decay. The key result proved by Gr\"ochenig in \cite{Grochenig_2006_Time} is a characterization of this type: a symbol belongs to the (weighted) Sj\"ostrand's class if and only if time-frequency shifts are almost eigenvectors of the corresponding Weyl operator. More precisely, the claim is the following.
 
 \begin{theorem}
 	Let $v$ be an admissible weight and fix a non-zero window $\varphi\in M^1_v(\mathbb{R}^{d})\left(\mathbb{R}^{d}\right)$ such that 
 	$\mathcal{G}\left(\varphi,\Lambda\right)$ is a Gabor frame  for $L^{2}\left(\mathbb{R}^{d}\right)$.
 	The following properties are equivalent:
 	\begin{enumerate}
 		\item[$(i)$] $\sigma\in M^{\infty,1}_{1\otimes v\circ J^{-1}}\left(\mathbb{R}^{2d}\right)$.
 		\item[$(ii)$] $\sigma\in\mathcal{S}'\left(\mathbb{R}^{2d}\right)$ and there exists
 		a function $H\in L^{1}_v\left(\mathbb{R}^{2d}\right)$ such
 		that
 		\[
 		\left|\left\langle \mathrm{Op_W}\left(\sigma\right)\pi\left(z\right)\varphi,\pi\left(w\right)\varphi\right\rangle \right|\le H\left(w-z\right),\qquad\forall w,z\in\mathbb{R}^{2d}.
 		\]
 		\item[$(iii)$] $\sigma\in\mathcal{S}'\left(\mathbb{R}^{2d}\right)$ and there exists
 		a sequence $h\in\ell^{1}_v\left(\Lambda\right)$ such that
 		\[
 		\left|\left\langle \mathrm{Op_W}\left(\sigma\right)\pi\left(\mu\right)\varphi,\pi\left(\lambda\right)\varphi\right\rangle \right|\le h\left(\lambda-\mu\right),\qquad\forall\lambda,\mu\in\Lambda.
 		\]
 	\end{enumerate}
 \end{theorem}
 
 This characterization is very strong: in particular, by applying Schwartz's kernel theorem, we also have:
 
 \begin{corollary}\label{weylchar}
 	Under the hypotheses of the previous Theorem, assume that \, $T:\mathcal{S}\left(\mathbb{R}^{d}\right)\rightarrow\mathcal{S}'\left(\mathbb{R}^{d}\right)$
 	is continuous and satisfies one of the following conditions:
 	\begin{itemize}
 		\item[$(i)$] $\left|\left\langle T\pi\left(z\right)\varphi,\pi\left(w\right)\varphi\right\rangle \right|\le H\left(w-z\right),\quad\forall w,z\in\mathbb{R}^{2d}$
 		for some $H\in L^{1}$. 
 		\item[$(ii)$] $\left|\left\langle T\pi\left(\mu\right)\varphi,\pi\left(\lambda\right)\varphi\right\rangle \right|\le h\left(\lambda-\mu\right),\quad\forall\lambda,\mu\in\Lambda$
 		for some $h\in\ell^{1}$. 
 	\end{itemize}
 	Therefore, $T=\mathrm{Op_W}\left(\sigma\right)$ for some symbol
 		$\sigma\in M^{\infty,1}_{1\otimes v\circ J^{-1}}\left(\mathbb{R}^{2d}\right).$
 \end{corollary}
 
 The proof of the main result heavily relies on a simple but crucial interplay between the entries of the channel matrix of $\mathrm{Op_W}$ and the short-time Fourier transform of the symbol, which will be discussed in complete generality in the subsequent section. 
 We mention that at this point Gr\"ochenig establishes a strong link with matrix algebra, hence heading towards a more conceptual discussion of the almost diagonalization property. In particular, it is easy to prove that $\sigma \in M^{\infty,1}_{1\otimes v\circ J^{-1}}$ if and only if its channel matrix $M(\sigma)$ belongs to the class $\mathcal{C}_v (\Lambda)$ of matrices $A=(a_{\lambda,\mu})_{\lambda,\mu\in\Lambda}$ such that there exists a sequence $h\in\ell^1_v$ which almost diagonalizes its entries, i.e. $$\| a_{\lambda,\mu}\| \leq h(\lambda-\mu),\qquad \lambda,\mu\in\Lambda.$$ It can be proved that $\mathcal{C}_v(\Lambda)$ is indeed a Banach *-algebra and this insight allows a natural extension if one considers other matrix algebras and investigates the relation between symbols and the membership of
their Gabor matrices in a matrix algebra. For further investigations in more general contexts, see for instance \cite{grochenig2008banach}.

Thanks to this fresh new formulation, the proofs of Sj\"ostrand's results provided by Gr\"ochenig are to certain extent more natural. Furthermore, they extend the previous ones since weighted spaces are considered. We summarize the main outcomes in the following claims. 

\begin{theorem}[Boundedness]~\\
	If $\sigma \in M^{\infty,1}_{1\otimes v\circ J^{-1}}$, then $\mathrm{Op_W}(\sigma)$ is bounded on $M^{p,q}_{m}$ for any $1\leq p,q \leq \infty$ and any $m\in\mathcal{M}_v$. In particular, if $\sigma \in M^{\infty,1}$, $\mathrm{Op_W}(\sigma)$ is bounded on $L^2(\mathbb{R}^{d})$ and
	\begin{itemize}
		\item if $1\leq p\leq2$, $\mathrm{Op_W}(\sigma)$ maps $L^p$ into $M^{p,p'}$;
		\item if $2\leq p \leq \infty$, $\mathrm{Op_W}(\sigma)$ maps $L^p$ into $M^p$. 
	\end{itemize}
\end{theorem}

\begin{theorem}[Algebra property]~\\
	If $v$ is a submultiplicative on $\mathbb{R}^{2d}$, then $M^{\infty,1}_v$ is a Banach $*$-algebra with respect to the twisted product $\sharp$ and the involution $\sigma \mapsto \overline{\sigma}$.
\end{theorem}

\begin{theorem}[Wiener property]~\\
	Assume that $v$ is a submultiplicative weight on $\mathbb{R}^{2d}$.  $\mathrm{Op_W}\left(M^{\infty,1}_v\right)$ is inverse-closed in $\mathcal{B}(L^2(\mathbb{R}^{d}))$ (i.e. if $\sigma \in \left(M^{\infty,1}_v\right)$ and $\mathrm{Op_W}(\sigma)$ is invertible on $L^2$, then $[\mathrm{Op_W}(\sigma)]^{-1}=\mathrm{Op_W}(\rho)$ for some $\rho\in\left(M^{\infty,1}_v\right)$) if and only if  $v$ satisfies the GRS condition \eqref{GRS}.
\end{theorem}

\begin{corollary}[Spectral invariance on modulation spaces]~\\
	Assume that $v$ is an admissible weight, $\sigma \in \left(M^{\infty,1}_v\right)$ and $\mathrm{Op_W}(\sigma)$ is invertible on $L^2$. Then, $\mathrm{Op_W}(\sigma)$ is simultaneously invertible on every modulation space $M^{p,q}_m(\mathbb{R}^{d})$, for any $1\leq p,q \leq \infty$ and $m\in\mathcal{M}_v$.
\end{corollary}

\begin{remark}
The intuition behind the last result is that the spectrum of an operator with suitably likable properties does not truly depend on the space on which it acts. 
In order to establish a link with Beals' theorem on spectral invariance in the context of classical pseudodifferential operators, notice that H\"ormander's class $$S^0_{0,0}(\mathbb{R}^{2d})=\{\sigma\in C^{\infty}(\mathbb{R}^{2d}) : \partial^{\alpha}\sigma \in L^{\infty}(\mathbb{R}^{2d}) \forall \alpha \in \mathbb{N}_0^{2d}\}$$ can be recast as intersection of Sj\"ostrand's classes with polynomial weights (cf. \cite{grochenig2008banach}), namely
$$  S^0_{0,0}(\mathbb{R}^{2d}) = \bigcap_{s\geq 0} M^{\infty,1}_{v_s}(\mathbb{R}^{2d}).$$
The Wiener property of these spaces leads to the conclusion that $\mathrm{Op_W}\left(S^0_{0,0}\right)$ is inverse-closed in $\mathcal{B} (L^2)$ too.
\end{remark}

\section{Almost diagonalization of $\tau$-pseudodifferential operators}
In a recent joint work of the author with E. Cordero and F. Nicola, an attempt has been made to follow the path outlined by Gr\"ochenig. The two directions investigated are 
\begin{enumerate}
	\item the extension of the almost-diagonalization theorem to more general operators;
	\item the search of an almost-diagonalization-like characterization of other symbol classes. 
\end{enumerate}
For what concerns the first point, $\tau$-pseudodifferential operators were investigated instead of those of Weyl type. We already discussed in the Introduction how this general class of operators extends in a natural way the previous one, which can be recovered as the case $\tau=1/2$. We were able to obtain an identical result with an identical proof - apart from the substantial modifications in the preliminary lemmas - see \cite{CNT18} for the details.

\begin{theorem}\label{almdiagsj}
	Let $v$ be an admissible weight on $\mathbb{R}^{2d}$. Consider  $\varphi\in M^1_v\left(\mathbb{R}^{d}\right)\setminus\{0\}$ and a lattice $\Lambda \subseteq \mathbb{R}^{2d}$ such that  $\mathcal{G}\left(\varphi,\Lambda\right)$ is a  Gabor frame for $L^{2}\left(\mathbb{R}^{d}\right)$.
	For any $\tau\in\left[0,1\right]$, the following properties are equivalent:
	\begin{enumerate}
		\item[$(i)$] $\sigma\in M_{1\otimes v\circ J^{-1}}^{\infty,1}\left(\mathbb{R}^{2d}\right)$.
		\item[$(ii)$] $\sigma\in\mathcal{S}'\left(\mathbb{R}^{2d}\right)$ and there exists
		a function $H_\tau \in L_{v}^{1}\left(\mathbb{R}^{2d}\right)$ such that
		$$
		\left|\left\langle \mathrm{Op}_{\tau} \left(\sigma\right)\pi\left(z\right)\varphi,\pi\left(w\right)\varphi\right\rangle \right|\le H_\tau\left(w-z\right)\qquad\forall w,z\in\mathbb{R}^{2d}.\label{eq:almdiag J}
		$$
		\item[$(iii)$] $\sigma\in\mathcal{S}'\left(\mathbb{R}^{2d}\right)$ and there exists
		a sequence $h_\tau \in\ell_{v}^{1}\left(\Lambda\right)$ such that
		$$
		\left|\left\langle \mathrm{Op}_{\tau} \left(\sigma\right)\pi\left(\mu\right)\varphi,\pi\left(\lambda\right)\varphi\right\rangle \right|\le h_\tau\left(\lambda-\mu\right)\qquad\forall\lambda,\mu\in\Lambda.\label{eq:almdiag discr}
		$$
	\end{enumerate}
\end{theorem}

This result is not surprising for at least two reasons. Looking at the mapping relating the symbols of different $\tau$-quantizations, namely (see for instance \cite{horm3,toft1})
$$ \operatorname*{Op}\nolimits_{\tau_1}(a_1)= \operatorname*{Op}\nolimits_{\tau_2}(a_2)\, \Leftrightarrow\,\widehat{a_2}(\xi_1,\xi_2)=e^{-2\pi i(\tau_2-\tau_1)\xi_1\xi_2}\widehat{a_1}(\xi_1,\xi_2),$$
we see that the map that relates a Weyl symbol to its $\tau$-counterpart is bounded in the Sj\"ostrand's class. At a more fundamental level, it is instructive to give a look at the crucial ingredient of the proof, which is the relation between the channel matrix of the $\tau$-pseudodifferential operator and the short-time Fourier transform of the symbol. 

\begin{proposition}
	\label{lem:STFT-gaborm}Fix a non-zero window $\varphi\in \mathcal{S}(\mathbb{R}^{d})$
	and set $\Phi_{\tau}=W_{\tau}\left(\varphi,\varphi\right)$ for $\tau\in\left[0,1\right]$.
	Then, for $\sigma\in \mathcal{S}'\left(\mathbb{R}^{2d}\right)$,
	\begin{equation}
	\left|\left\langle \mathrm{Op}_{\tau} \left(\sigma\right)\pi\left(z\right)\varphi,\pi\left(w\right)\varphi\right\rangle \right|=\left|{V}_{\Phi_{\tau}}\sigma\left(\mathcal{T}_{\tau}\left(z,w\right),J\left(w-z\right)\right)\right|=\left|{V}_{\Phi_{\tau}}\sigma\left(x,y\right)\right|\label{eq:gaborm as STFT}
	\end{equation}
	and
	\begin{equation}
	\left|{V}_{\Phi_{\tau}}\sigma\left(x,y\right)\right|=\left|\left\langle \mathrm{Op}_{\tau} \left(\sigma\right)\pi\left(z\left(x,y\right)\right)\varphi,\pi\left(w\left(x,y\right)\right)\varphi\right\rangle \right|,\label{eq:STFT as gaborm}
	\end{equation}
	for all $w,z,x,y\in\mathbb{R}^{2d}$, where $\mathcal{T}_\tau$ is defined as \begin{equation}\label{Ttau}
	\mathcal{T}_{\tau}\left(z,w\right)=\left(\begin{array}{c}
	\left(1-\tau\right)z_{1}+\tau w_{1}\\
	\tau z_{2}+\left(1-\tau\right)w_{2}
	\end{array}\right)\quad z=(z_1,z_2),\,w=(w_1,w_2)\in\mathbb{R}^{2d}.
	\end{equation} and
	\begin{equation}
	z\left(x,y\right)=\left(\begin{array}{c}
	x_{1}+\left(1-\tau\right)y_{2}\\
	x_{2}-\tau y_{1}
	\end{array}\right),\qquad w\left(x,y\right)=\left(\begin{array}{c}
	x_{1}-\tau y_{2}\\
	x_{2}+\left(1-\tau\right)y_{1}
	\end{array}\right).\label{eq:z,w}
	\end{equation}
\end{proposition}

The main remark here is that the controlling function $H_{\tau}\in L^1_v(\mathbb{R}^{d})$ in the almost diagonalization theorem can be chosen as the so-called \emph{grand symbol} associated to $\sigma\in M^{\infty,1}_{v\circ J^{-1}}$ (according to \cite{grocomp}): for the general $\tau$-case, we have 
\[
H_{\tau}(v)=\sup_{u\in\mathbb{R}^{2d}} \left| V_{\Phi_{\tau}}\sigma(u,Jv)\right|.
\]
The choice of the grand symbol is quite natural if one looks at the modulation norm in the Sj\"ostrand's class. However, it is clear that the dependence from $\tau$ is completely confined to the window function $\Phi_\tau$ and does not affect the variable $v\in\mathbb{R}^{2d}$, which corresponds to the frequency variable for the short-time Fourier transform of the symbol. The proof of the general case can thus proceed exactly as the one for Weyl case. 
We remark that also Corollary \ref{weylchar} generalizes in the obvious way. 

It is reasonable at this stage to ask what happens if a slight modification of the grand symbol is taken into account, that is: what happens if we look at the time dependence of $V_{\Phi_{\tau}} \sigma$? This is equivalent to wonder if similar arguments extend in some fashion to Fourier transform of symbols in the Sj\"ostrand's class, namely symbols in a suitably weighted version of Wiener amalgam space $W\left(\mathcal{F} L^{\infty},L^1\right)=\mathcal{F} M^{\infty,1}$ - hereinafter referred to $\mathcal{F}$-Sj\"ostrand's class. The main outcome we got is the following.

\begin{theorem}\label{thm:almost diag}
	Let $v$ be an admissible weight function on $\mathbb{R}^{2d}$. Consider  $\varphi\in M^1_v\left(\mathbb{R}^{d}\right)\setminus\{0\}$.
	For any $\tau\in\left(0,1\right)$, the following properties are equivalent:
	\begin{enumerate}
		\item[$(i)$] $\sigma\in W\left(\mathcal{F}L^{\infty},L_{v\circ\mathcal{B}_{\tau}}^{1}\right)\left(\mathbb{R}^{2d}\right)$.
		\item[$(ii)$] $\sigma\in\mathcal{S}'\left(\mathbb{R}^{2d}\right)$ and there exists
		a function $H_\tau\in L_{v}^{1}\left(\mathbb{R}^{2d}\right)$ such that
		\begin{equation}\label{eq:almdiagfsj}
		\left|\left\langle \mathrm{Op}_{\tau} \left(\sigma\right)\pi\left(z\right)\varphi,\pi\left(w\right)\varphi\right\rangle \right|\le H_\tau\left(w-\mathcal{U}_{\tau}z\right)\qquad\forall w,z\in\mathbb{R}^{2d},
		\end{equation}
		where the matrices $\mathcal{B}_{\tau} $ and  $\mathcal{U}_{\tau}$ are defined as 
		\begin{equation}\label{btau utau}
		\mathcal{B}_{\tau}=\left(\begin{array}{cc}
		\frac{1}{1-\tau}I_{d\times d} & 0_{d\times d}\\
		0_{d\times d} & \frac{1}{\tau}I_{d\times d}
		\end{array}\right), \qquad \mathcal{U}_{\tau}=-\left(\begin{array}{cc}
		\frac{\tau}{1-\tau}I_{d\times d} & 0_{d\times d}\\
		0_{d\times d} & \frac{1-\tau}{\tau}I_{d\times d}
		\end{array}\right)\in\mathrm{Sp}\left(2d,\mathbb{R}\right).
		\end{equation}
If $\tau\in [0,1]$, the estimate in \eqref{eq:almdiagfsj} weakens as follows:
		\begin{enumerate}
					\item[$(ii')$] $\sigma\in\mathcal{S}'\left(\mathbb{R}^{2d}\right)$ and there exists
			a function $H_\tau\in L_{v}^{1}\left(\mathbb{R}^{2d}\right)$ such that
			\begin{equation}\label{weakalmd}
			\left|\left\langle \mathrm{Op}_{\tau} \left(\sigma\right)\pi\left(z\right)\varphi,\pi\left(w\right)\varphi\right\rangle \right|\leq H_\tau\left(\mathcal{T}_{\tau}(w,z)\right)\qquad\forall w,z\in\mathbb{R}^{2d}.
			\end{equation}
		\end{enumerate}  
	\end{enumerate}
\end{theorem}

A number of differences arise with respect to its counterpart for Sj\"ostrand's symbols. First, the almost diagonalization of the (continuous) channel matrix is lost, but this is still a well-organized matrix: in the favourable case $\tau=(0,1)$,  \eqref{eq:almdiagfsj} can be interpreted as a measure of the concentration of the time-frequency representation of $\mathrm{Op}_{\tau} (\sigma)$ along the graph of the map $\mathcal{U}_{\tau}$. If we include the endpoints, the estimate loses this meaning too.

Furthermore, notice that the discrete characterization via Gabor frames is lost, the main obstruction being the following: for a given lattice $\Lambda$, the inclusion $\mathcal{U}_\tau \Lambda \subseteq \Lambda$ holds if and only if  $\tau=1/2$, i.e. $\mathcal{U}_\tau=\mathcal{U}_{1/2}=-I_{2d\times 2d}$.  In this particular framework, the matrix $\mathcal{B}_{1/2}$ then becomes
$\mathcal{B}_{1/2}=2 I_{2d\times 2d}$ and the symmetry of Weyl operators is rewarded by an additional characterization:
\begin{enumerate}
	\item[$(iii')$] $\sigma\in\mathcal{S}'\left(\mathbb{R}^{2d}\right)$ and there exists
	a sequence $h \in\ell_{v}^{1}\left(\Lambda\right)$ such that
	$$
	\left|\left\langle \mathrm{Op_W}\left(\sigma\right)\pi\left(\mu\right)\varphi,\pi\left(\lambda\right)\varphi\right\rangle \right|\le h \left(\lambda+\mu\right)\qquad\forall\lambda,\mu\in\Lambda.\label{eq:almdiag discr2}
	$$
\end{enumerate}

\section{Consequences of almost diagonalization}
\subsection{Boundedness}
We are now able to study the boundedness of $\tau$-pseudodifferential operators covering several possible choices for symbols classes and spaces on which they act. 
If one considers the action of $\tau$-pseudodifferential operators on modulation spaces, a Sj\"ostrand-type result for symbols in the Sj\"ostrand's class can be inferred by means of the same arguments applied in the Weyl case. 
\begin{theorem}\label{TW}
	Consider  $m\in\mathcal{M}_{v}\left(\mathbb{R}^{2d}\right)$ satisfying \eqref{M}. For any $\tau\in [0,1]$ and $\sigma\in M^{\infty,1}_{1\otimes v\circ J^{-1}}$ the operator $\mathrm{Op}_{\tau} (\sigma)$ is bounded on $M^{p,q}_m(\mathbb{R}^{d})$, and there exists a constant $C_\tau>0$ such that
	\begin{equation}\label{Toftweight}
	\| \mathrm{Op}_{\tau} (\sigma)\|_{M^{p,q}_m}\leq C_\tau \|\sigma\|_{M^{\infty,1}_{1\otimes v\circ J^{-1}}}.
	\end{equation}
\end{theorem}

In order to address the problem of boundedness of $\tau$-pseudodifferential operators on modulation spaces with symbols in $\mathcal{F}$-Sj\"ostrand's class, a different strategy is needed. Following \cite{generalizedmetaplectic}, the idea is to recast $\mathrm{Op}_{\tau} (\sigma)$ as the transformation (via the short-time Fourier transform and its adjoint) of an integral operator with the channel matrix as distributional kernel. Therefore, the almost diagonalization property allows to obtain the desired estimates and claim the following result. 

\begin{theorem}
	\label{thm:bounded optau} Fix $m\in\mathcal{M}_{v}$ satisfying \eqref{M}. For $\tau\in\left(0,1\right)$
	consider a symbol $\sigma\in W(\mathcal{F}L^{\infty},L_{v\circ\mathcal{B}_{\tau}}^{1})\left(\mathbb{R}^{2d}\right)$, with the matrix $\mathcal{B}_{\tau}$ defined in \eqref{btau utau}. Then the operator $\mathrm{Op}_{\tau} (\sigma)$ is bounded from $M_{m}^{p,q}\left(\mathbb{R}^{d}\right)$
	to $M_{m\circ\mathcal{U}_{1-\tau}^{-1}}^{p,q}\left(\mathbb{R}^{d}\right)$,
	$1\le p,q\le\infty$. 
\end{theorem}

We now turn to consider the boundedness of $\tau$-pseudodifferential operators on Wiener amalgam spaces. Looking for a big picture and given that modulation and Wiener amalgam spaces are intertwined by the Fourier transform, it is natural to wonder if continuity properties of an operator acting on modulation spaces may still hold true when it acts on the corresponding amalgam spaces. In the case of $\tau$-pseudodifferential operators the answer is yes but heavily relies on the particular way Fourier transform and $\tau$-pseudodifferential operators commute. This phenomenon is a special case of the symplectic covariance property of Shubin calculus, which we briefly recall - see \cite{dgTRAN2013} for a comprehensive discussion on the issue.

\begin{lemma}
	\label{lem:symplectic covariance tau}For any $\sigma\in\mathcal{S}'\left(\mathbb{R}^{2d}\right)$
	and $\tau\in\left[0,1\right]$, 
	\[
	\mathcal{F}\mathrm{Op}_{\tau} \left(\sigma\right)\mathcal{F}^{-1}=\mathrm{Op}_{1-\tau}\left(\sigma\circ J^{-1}\right).
	\]
\end{lemma}

This property, along with other preliminary results, allows to quickly prove the desired claims for symbols in both Sj\"ostrand's class and the corresponding amalgam space. 

\begin{theorem}
	Consider  $m=m_{1}\otimes m_{2}\in\mathcal{M}_{v}\left(\mathbb{R}^{2d}\right)$ satisfying \eqref{M}.
	For any $\tau\in [0,1]$ and $\sigma\in M_{1\otimes v }^{\infty,1}\left(\mathbb{R}^{2d}\right)$,
	the operator $\mathrm{Op}_{\tau} (\sigma)$ is bounded on $W\left(\mathcal{F}L_{m_{1}}^{p},L_{m_{2}}^{q}\right)\left(\mathbb{R}^{d}\right)$
	with
	$$
	\| \mathrm{Op}_{\tau} (\sigma)\|_{W\left(\mathcal{F}L_{m_{1}}^{p},L_{m_{2}}^{q}\right)}\leq C_\tau \|\sigma\|_{M_{1\otimes v}^{\infty,1}},$$
	for a suitable $C_\tau>0$.
\end{theorem}

\begin{theorem}
	Consider  $m=m_{1}\otimes m_{2}\in\mathcal{M}_{v}\left(\mathbb{R}^{2d}\right)$ satisfying \eqref{M}. For any $\tau\in\left(0,1\right)$ and $\sigma\in W\left(\mathcal{F}L^{\infty},L_{v\circ\mathcal{B}_{\tau}\circ J^{-1}}^{1}\right)\left(\mathbb{R}^{2d}\right)$, the operator $\mathrm{Op}_{\tau}\sigma$ is bounded from $W\left(\mathcal{F}L_{m_{1}}^{p},L_{m_{2}}^{q}\right)\left(\mathbb{R}^{d}\right)$ to $W\left(\mathcal{F}L_{m_{1}\circ\left(\mathcal{U}_{1-\tau}^{-1}\right)_{1}}^{p},L_{m_{2}\circ\left(\mathcal{U}_{1-\tau}^{-1}\right)_{2}}^{q}\right)\left(\mathbb{R}^{d}\right)$, $1\le p,q\le\infty$, where 
	\[ \left(\mathcal{U}_{1-\tau}^{-1}\right)_{1}\left(x\right)=-\frac{\tau}{1-\tau}x,\qquad\left(\mathcal{U}_{1-\tau}^{-1}\right)_{2}\left(x\right)=-\frac{1-\tau}{\tau}x,\qquad x\in\mathbb{R}^{d}.
	\]
\end{theorem} 

We finally remark that even if the results with symbols in $\mathcal{F}$-Sj\"ostrand's class do not hold for the endpoint cases $\tau=0$ and $\tau=1$, it is still possible to use the weak characterization \eqref{weakalmd} to construct ad hoc examples of bounded operators. 

\begin{proposition}
	Assume $\sigma \in W(\mathcal{F} L^\infty, L^1)(\mathbb{R}^{2d})$.  
	\begin{enumerate}
		\item The  Kohn-Nirenberg operator $\mathrm{Op_{KN}}(\sigma) \,\, (\tau = 0)$ is bounded on $M^{1,\infty}(\mathbb{R}^{d})$.
		\item The anti-Kohn-Nirenberg $\mathrm{Op}_{1}(\sigma) \,\, (\tau = 1)$ is bounded on $W(\mathcal{F} L^1, L^\infty)(\mathbb{R}^{d})$. 
	\end{enumerate}
	
\end{proposition}

\subsection{Algebra and Wiener properties}

To conclude, we give a brief summary on the extension of the other properties studied by Sj\"ostrand, namely algebra and Wiener property, to $\tau$-pseudodifferential operators. 
Wiener algebras of pseudodifferential operators have been already investigated by Cordero, Gr\"ochenig, Nicola and Rodino in several occasions, see for instance \cite{cordero2013wiener,generalizedmetaplectic}. Let us recall the definition and the relevant properties of generalized metaplectic operators, introduced by the aforementioned authors.

\begin{definition}
	Given $\mathcal{A} \in \mathrm{Sp}\left(2d,\mathbb{R}\right)$,
	$g\in\mathcal{S}(\mathbb{R}^{d})$,  and $s\geq0$, a
	linear operator $T:\mathcal{S}(\mathbb{R}^{d})\to\mathcal{S}'(\mathbb{R}^{d})$ belongs to the
	class $FIO(\mathcal{A},v_s)$ of generalized metaplectic operators if $$\exists H\in L^1_{v_s}(\mathbb{R}^{2d}) \text{ such that }
	|\langle T \pi(z) g,\pi(w)g\rangle|\leq H(w-\mathcal{A} z),\qquad \forall w,z\in\mathbb{R}^{2d}.
	$$
	
\end{definition}

\begin{theorem}
	
	Fix $\mathcal{A}_{i}\in\mathrm{Sp}\left(2d,\mathbb{R}\right)$, $s_{i}\ge0$,
	$m_{i}\in\mathcal{M}_{v_{s_{i}}}$, and $T_{i}\in FIO\left(\mathcal{A}_{i},v_{s_{i}}\right)$,
	$i=0,1,2$. 
	\begin{enumerate}
		\item $T_{0}$ is bounded from $M_{m_{0}}^{p}\left(\mathbb{R}^{d}\right)$
		to $M_{m_{0}\circ\mathcal{A}_{i}^{-1}}^{p}\left(\mathbb{R}^{d}\right)$
		for any $1\le p\le\infty$.
		\item $T_{1}T_{2}\in FIO\left(\mathcal{A}_{1}\mathcal{A}_{2},v_{s}\right)$,
		where $s=\min\left\{ s_{1},s_{2}\right\} $.
		\item If $T_{0}$ is invertible in $L^{2}\left(\mathbb{R}^{d}\right)$,
		then $T_{0}^{-1}\in FIO\left(\mathcal{A}_{0}^{-1},v_{s_{0}}\right)$.
	\end{enumerate}
	
\end{theorem}

In short, the class 
$$FIO(\mathrm{Sp}(2d,\mathbb{R}),v_s)=\bigcup_{\mathcal{A}\in \mathrm{Sp}(2d,\mathbb{R})} FIO(\mathcal{A},v_s)$$ is a Wiener sub-algebra of $\mathcal{B}(L^2(\mathbb{R}^{d}))$.
	In view of the defining property of operators in $FIO(\mathcal{A},v_s)$, we immediately recognize that for any $\tau \in (0,1)$, if $\sigma \in W(\mathcal{F} L^\infty, L^1_{v_s})$ then $\mathrm{Op}_\tau(\sigma) \in FIO(\mathcal{U}_{\tau},v_s)$. Therefore, if we limit to consider admissible weights of polynomial type $v_s$ on $\mathbb{R}^{d}$, $s\ge0$, we are able to establish a fruitful connection and to derive a number of properties without any effort. For instance, we have another boundedness result.
\begin{corollary}
	If  $\sigma \in W(\mathcal{F} L^\infty, L^1_{v_s})(\mathbb{R}^{2d})$, $s\geq 0$, then the operator $\mathrm{Op}_{\tau} (\sigma)$ is bounded on every modulation space $M^p_{v_s}(\mathbb{R}^{d})$, for $1\leq p\leq \infty$ and $\tau \in (0,1)$.
\end{corollary}

For what concerns the algebra property, we in fact have a no-go result. By inspecting the composition properties of matrices $\mathcal{U}_{\tau}$, we notice that there is no $\tau\in\left(0,1\right)$ such that $\mathcal{U}_{\tau_{1}}\mathcal{U}_{\tau_{2}}=\mathcal{U}_{\tau}$.
This implies that there is no $\tau$-quantization rule
such that composition of $\tau$-operators with symbols in $W\left(\mathcal{F}L^{\infty},L_{v_{s}}^{1}\right)$
has symbol in the same class. We can only state weaker algebraic results, such as the following property of ``symmetry" with respect to the Weyl quantization. 
\begin{theorem}
	For any $a,b\in W\left(\mathcal{F}L^{\infty},L_{v_{s}}^{1}\right)(\mathbb{R}^{2d})$
	and $\tau\in\left(0,1\right)$, there exists a symbol $c\in M_{1\otimes v_{s}}^{\infty,1}(\mathbb{R}^{2d})$
	such that
	\[
	\mathrm{Op}_{\tau}\left(a\right)\mathrm{Op}_{1-\tau}\left(b\right)=\mathrm{Op}_{1/2}\left(c\right).
	\]
\end{theorem} 

Also notice that, given $a\in W\left(\mathcal{F}L^{\infty},L_{v_{s}}^{1}\right)$,
$b\in M_{1\otimes v_{s}}^{\infty,1}$ and $\tau,\tau_{0}\in\left(0,1\right)$,
we have 
\[
\mathrm{Op}_{\tau_{0}}\left(b\right)\mathrm{Op}_{\tau}\left(a\right)=\mathrm{Op}_{\tau}\left(c_{1}\right),\qquad\mathrm{Op}_{\tau}\left(a\right)\mathrm{Op}_{\tau_{0}}\left(b\right)=\mathrm{Op}_{\tau}\left(c_{2}\right),
\]
for some $c_{1},c_{2}\in W\left(\mathcal{F}L^{\infty},L_{v_{s}}^{1}\right)$.
This means that, for fixed quantization rules $\tau,\tau_{0}$, the
amalgam space $W\left(\mathcal{F}L^{\infty},L_{v_{s}}^{1}\right)(\mathbb{R}^{2d})$
is a bimodule over the algebra $M_{1\otimes v_{s}}^{\infty,1}(\mathbb{R}^{2d})$ under
the laws 
\[
M_{1\otimes v_{s}}^{\infty,1}\times W\left(\mathcal{F}L^{\infty},L_{v_{s}}^{1}\right)\rightarrow W\left(\mathcal{F}L^{\infty},L_{v_{s}}^{1}\right)\,\,:\,\,\left(b,a\right)\mapsto c_{1},
\]
\[
W\left(\mathcal{F}L^{\infty},L_{v_{s}}^{1}\right)\times M_{1\otimes v_{s}}^{\infty,1}\rightarrow W\left(\mathcal{F}L^{\infty},L_{v_{s}}^{1}\right)\,\,:\,\,\left(a,b\right)\mapsto c_{2},
\]
with $c_{1}$ and $c_{2}$ as before. 

Finally, after noticing that $\mathcal{U}_{\tau}^{-1}=\mathcal{U}_{1-\tau}$
for any $\tau\in\left(0,1\right)$, a Wiener-like property comes at the price of passing to the complementary $\tau$-quantization when inverting $\mathrm{Op}_{\tau}$. 

\begin{theorem}
	For any $\tau\in\left(0,1\right)$ and $a\in W\left(\mathcal{F}L^{\infty},L_{v_{s}}^{1}\right)(\mathbb{R}^{2d})$
	such that $\mathrm{Op}_{\tau}\left(a\right)$ is invertible on $L^{2}\left(\mathbb{R}^{d}\right)$,
	we have 
	\[
	\mathrm{Op}_{\tau}\left(a\right)^{-1}=\mathrm{Op}_{1-\tau}\left(b\right)
	\]
	for some $b\in W\left(\mathcal{F}L^{\infty},L_{v_{s}}^{1}\right)(\mathbb{R}^{2d})$.
\end{theorem}

\end{document}